\documentclass[12pt]{amsart}

\usepackage{graphicx}
\usepackage{amsthm}
\usepackage{amsmath, amssymb}
\swapnumbers

\theoremstyle{plain}
\newtheorem{Thm}{Theorem}

\newtheorem{Lem}[Thm]{Lemma}
\newtheorem*{Claim}{Claim}

\theoremstyle{definition}

\begin{document}

\title{Horizontal Heegaard splittings of Seifert fibered spaces}

\author{Jesse Johnson}
\address{\hskip-\parindent
        Department of Mathematics \\
        Yale University \\
        PO Box 208283 \\
        New Haven, CT 06520 \\
        USA}
\email{jessee.johnson@yale.edu}


\thanks{Research supported by NSF MSPRF grant 0602368}

\begin{abstract}
We show that if an orientable Seifert fibered space $M$ with an orientable genus $g$ base space admits a strongly irreducible horizontal Heegaard splitting then there is a one-to-one correspondence between isotopy classes of strongly irreducible horizontal Heegaard splittings and elements of $\mathbf{Z}^{2g}$.   The correspondence is determined by the slopes of intersection of each Heegaard splitting with a collection of $2g$ incompressible tori in $M$.  We also show that there are Seifert fibered spaces with infinitely many non-isotopic Heegaard splittings that determine Nielsen equivalent generating systems for the fundamental group of $M$.
\end{abstract}

\maketitle

\section{Introduction}
\label{introsect}
Certain closed Seifert Fibered spaces are known to admit a type of Heegaard splitting called a horizontal Heegaard splitting.  Bachman and Derby-Talbot~\cite{bdt} showed that any Seifert fibered that admits a strongly irreducible horizontal splitting admits infinitely many isotopy classes of horizontal splittings.  We improve their analysis to show the following:

Let $M$ be an orientable Seifert fibered space with  base space an orientable genus $g$ surface and let $T_1,\dots,T_{2g}$ be vertical tori in $M$ such that $T_{i} \cap T_{j}$ is a single loop for $i$ odd, $j = i+1$ (or vice versa) and empty otherwise.  The complement in $M$ of a regular neighborhood of these tori is a Seifert fibered space over a $g$-times punctured sphere.  

\begin{Thm}
\label{mainthm}
If $M$ admits a strongly irreducible horizontal Heegaard splitting then for every $2g$-tuple of integers $(s_1,\dots,s_{2g}) \in \mathbf{Z}^{2g}$, there is a unique (up to isotopy) strongly irreducible, horizontal Heegaard splitting that intersects each $T_i$ in a family of essential loops with slope $s_i$.  Moreover, Heegaard splittings that define distinct $2g$-tuples of slopes are not isotopic.
\end{Thm}

A \textit{Heegaard splitting} for a compact, closed, orientable 3-manifold $M$ is a triple $(\Sigma, H_1, H_2)$ where $\Sigma \subset M$ is a compact, closed, two-sided surface and $H_1, H_2 \subset M$ are handlebodies (homeomorphic copies of closed regular neighborhoods of connected, finite graphs in $S^3$) with $\partial H_1 = \Sigma = \partial H_2$ and $H_1 \cup H_2 = M$.  

A Heegaard splitting $(\Sigma, H_1, H_2)$ is \textit{strongly irreducible} if every essential, properly embedded disk in $H_1$ intersects every essential, properly embedded disk in $H_2$.  We will describe the construction of a horizontal Heegaard splitting in Section~\ref{sfsect}.

Given a Heegaard splitting $(\Sigma, H_1, H_2)$ of $M$, there is a smooth function $f : M \rightarrow [0,1]$ such that the pre-image of each point in $(0,1)$ is a surface isotopic to $\Sigma$ and the pre-images of $\{0\}$ and $\{1\}$ are graphs (called \textit{spines}) in $H_1$ and $H_2$.  Such a function is called a \textit{sweep-out} (see~\cite{mine}) and the restriction of $f$ to a vertical torus in $M$ is (generically) a Morse function.  A Morse function on a torus always has level sets that are essential in the torus.  Level sets of a Morse function are pairwise disjoint and disjoint essential loops in a torus are parallel so $f$ determines a unique isotopy class of simple closed curves in the torus.  

We will describe below how a simple closed curve in a vertical torus determines a rational number called its slope.  Different sweep-outs will restrict to different Morse functions on $T$, so a Heegaard splitting may determine more than one slope.  We will show that in many cases if two sweep-outs come from the same Heegaard splitting then they will determine the same slope on the vertical torus.  In particular, for $M$ a Seifert fibered space with orientable base space and $T_1,\dots,T_{2g}$ vertical tori in $M$ as above, we show the following:

\begin{Lem}
\label{oneslopelem}
A strongly irreducible Heegaard splitting of a Seifert fibered space $M$ determines a unique slope in each vertical torus $T_i$.
\end{Lem}

This is proved in Section~\ref{sfsect}, based on techniques developed in Section~\ref{sumsect}, and shows one direction of Theorem~\ref{mainthm}.  The other direction follows from the construction of horizontal Heegaard splittings and is also proved in Section~\ref{sfsect}.

In Sections~\ref{dpsect} and~\ref{nielsect}, we consider the generating set for the fundamental group of $M$.  Two generating sets are called \textit{Nielsen equivalent} if one can be changed to the other by a finite number of type-one Tietze moves (i.e. replacing the $i$th generator with its inverse or with the product of the $i$th and the $j$th generator for some $i \neq j$).  

The fundamental group of each handlebody in a Heegaard splitting is a free group and the inclusion of its fundamental group into $\pi_1(M)$ determines a generating set for $\pi_1(M)$.  If the handlebodies of two Heegaard splittings determine generating sets for $\pi_1(M)$ that are not Nielsen equivalent then the Heegaard splittings can not be isotopic.  Lustig and Moriah have used Nielsen equivalence to distinguish vertical Heegaard splittings of Seifert fibered spaces~\cite{lm1} as well as Heegaard splittings of certain hyperbolic 3-manifolds~\cite{lm2}.  We show that unfortunately, Nielsen class does not always distinguish non-isotopic Heegaard splittings.  In particular, we describe in Section~\ref{nielsect} a family of Seifert fibered space over the torus with two singular fibers such that each admits infinitely many non-isotopic Heegaard splittings whose handlebodies determine Nielsen equivalent generating sets in $\pi_1(M)$.

\section{Toroidal summands}
\label{sumsect}

Let $M$ be a compact, closed, orientable, irreducible 3-manifold (not necessaily a Seifer fibered space) and let $N \subset M$ be a submanifold homeomorphic to $T \times S^1$ where $T$ is a once punctured torus.  Assume $\partial N$ is incompressible in $M$.  (If $\partial N$ is compressible in $M$ then it compresses to a sphere in the complement of $N$ so, because $M$ is irreducible, $N$ must be a solid torus.)

There are two canonical simple closed curves in $\partial N$ that are picked out by the topology: a meridian $\mu$ that is the boundary of an incompressible torus $T \times \{y\}$ in $T \times S^1$ (for some $y \in S^1$) and a longitude $\lambda$ that is the slope of a vertical loop $\{x\} \times S^1$ ($y \in \partial T$).  The meridian $\mu$ is the unique (up to isotopy) loop in $\partial N$ that is homology trivial in $N$, so it is determined independently of the product structure on $N$.  Every essential annulus properly embedded in $N$ has boundary parallel to $\lambda$ so this loop is also independent of the product structure.  Any simple closed curve in $\partial N$ is a sum $p \mu + q \lambda$ and thus determines a fraction $\frac{p}{q} \in \mathbf{Q} \cup \{\frac{1}{0}\}$, called its \textit{slope}.  

For any essential, simple closed curve $\ell$ in $T$, the subset $\ell \times S^1 \subset T \times S^1$ is a non-separating incompressible torus in $M$.  We can define slopes $\mu' = \ell \times \{y\}$ for $y \in S^1$ and $\lambda' = \{x\} \times S^1$ for $x \in \ell$, so again each loop in $\ell \times S^1$ determines a slope $\frac{p}{q}$.  In this case, the loop $\mu'$ is determined by the product structure of $N$, not the topology alone.  A different product structure will imply a different $\mu'$.  For our purposes, it suffices to fix a product structure on $N$, since we will always be dealing with these slopes in a relative way.

Let $(\Sigma, H_1, H_2)$ be a Heegaard splitting for $M$.  Let $f : M \rightarrow [0,1]$ be a sweep-out such that each level surface of $f$ is isotopic to $\Sigma$.  Let $S = \ell \times S^1$ be a vertical torus in $N$.  After an arbitrarily small isotopy of $f$, the restriction of $f$ to $S$ will be a Morse function.  As mentioned above, a Morse function on a torus always has an essential level set and the essential levels define a single isotopy class of simple closed curves.  Thus $f$ determines a unique slope in $S$.  

We will say that $\Sigma$ \textit{determines} a slope $\frac{p}{q}$ on $S$ if there is a sweep-out $f$ with level sets isotopic to $\Sigma$ such that the restriction of $f$ to $S$ has a level set in $S$ with slope $\frac{p}{q}$.  As noted above, a Heegaard splitting may determine more than one slope.  If the intersection of $\Sigma$ with $S$ contains an essential loop of slope $\frac{p}{q}$ then $f$ can be chosen so that $\Sigma$ is a level set of $f$ (rather than just isotopic to one) so $\Sigma$ determines the slope $\frac{p}{q}$.  Conversely, for any sweep-out for $\Sigma$, each level surface is isotopic to $\Sigma$.  Thus $\Sigma$ determines a slope $\frac{p}{q}$ in $S$ if and only if $\Sigma$ can be isotoped so that the intersection contains a loop with that slope.

\begin{Lem}
\label{slopechangelem}
If a strongly irreducible genus $g$ Heegaard splitting $(\Sigma, H_1, H_2)$ for $M$ determines more than one slope in a vertical torus $S$ in $N$ then $\Sigma$ can be isotoped so that the closure of $\Sigma \setminus N$ in the closure of $M \setminus N$ is a properly embedded incompressible genus $g - 3$ surface whose boundary is a pair of loops in $\partial N$, each with slope $1$ or $-1$.
\end{Lem}

Before we begin the proof, recall that a smooth function $f$ is \textit{Morse} if every critical point is non-degenerate and no two critical points are in the same level.  A function is \textit{near-Morse} if either all but one of its critical points are non-degenerate and all are in distinct levels, or all its critical points are non-degenerate and all but two are in distinct levels.

\begin{proof}
If $(\Sigma, H_1, H_2)$ determines more than one slope in $S$ then there are sweep-outs $f$, $f'$ such that $\Sigma$ is isotopic to both a level surface of $f$ and of $f'$ and such that the essential level sets of $f |_S$ and $f'|_S$ determine different slopes.  Because $f$ and $f'$ are sweep-outs for the same Heegaard splitting, there is an isotopy of $M$ taking a level surface of $f'$ to a level surface of $f$.  In particular, there is a family of sweep-outs $\{f_t | t \in [0,1]\}$ such that $f_0$ determines the same slope in $S$ as $f'$ and $f_1$ determines the same slope in $S$ as $f$.  

Assume the family of sweep-outs is generic with respect to $S$, i.e. that $f_t|_S$ is near-Morse for finitely many values of $t$ and Morse for the remaining values of $t$.  For any value $t_0$ such that $f_{t_0}|_S$ is a Morse function, there is a neighborhood of $t_0$ in $[0,1]$ such that for any $t$ in this neighborhood, $f_t|_S$ is isotopic (in $S$) to $f_{t_0}|_S$.  Thus the slope of the essential levels can only change at the near-Morse values of $t$.

If two essential loops in a torus are disjoint then they are parallel, and thus define the same slope.  Thus if the essential slope changes at a near-Morse value $t_0$ then the regular levels of $f_{t_0}|_S$ must all be trivial in $S$.  This is the case if and only if each component of the complement of the critical levels is contained in an open disk in $S$.  If $f_{t_0}|_S$ is a near-Morse function with a degenerate critical point (but its critical points are in distinct levels) then the complement of the critical levels must still contain a component not contained in a disk.  Thus the slope of the essential levels of $f_t|_S$ can only change when $f_{t_0}|_S$ has two saddles at the same level, and this level set cuts $S$ into disks.  

The critical level containing the two saddle singularities is a graph with two valence four vertices and thus four edges.  There are exactly two (homeomorphism classes of) connected graphs with four edges and two valence four vertices:  Let $\Gamma_0$ be a two-vertex graph in which two edges pass between the two vertices and one edge goes from each vertex back to itself.  Let $\Gamma_1$ be a two-vertex graph in which each edge goes from one vertex to the other.

Let $\Gamma$ be a critical level set of a near-Morse function on an oriented surface $S$ such that $\Gamma$ is homeomorphic to $\Gamma_0$ or $\Gamma_1$.  Given an orientation for an edge of $\Gamma$, the orientation of $S$ defines a transverse orientation.  Choose an orientation for each edge so that the transverse orientation points in the direction in which the near-Morse function is increasing.  The embedding of $\Gamma$ suggests a cyclic ordering of the ends of the edges that enter each vertex.  Because each vertex is at a saddle singularity, the edges must alternate whether they point towards the vertex or away.

If $\Gamma$ is homeomorphic to $\Gamma_0$ then for each edge that passes from a vertex to itself, one end points towards the vertex and the other away.  Thus the ends of each such edge are adjacent in the cyclic ordering around the vertex.  This implies that a regular neighborhood of $\Gamma$ is planar.  If $S$ is a torus then the complement of $\Gamma$ must contain a component that is not contained in a disk in $S$.

Thus if the slope defined by the Morse function changes, the level containing two saddles must be homeomorphic to $\Gamma_1$.  There is a unique (up to homeomorphism) way that such a graph can be embedded in a torus such that its complement is a collection of diks.  This is shown at the bottom left of Figure~\ref{torusfig}.  The top left picture shows the intersection of this level set with a square whose sides are glued to construct the torus, chosen so that the two vertices of $\Gamma$ are in the two edges of the glued square.
\begin{figure}[htb]
  \begin{center}
  \includegraphics[width=3.5in]{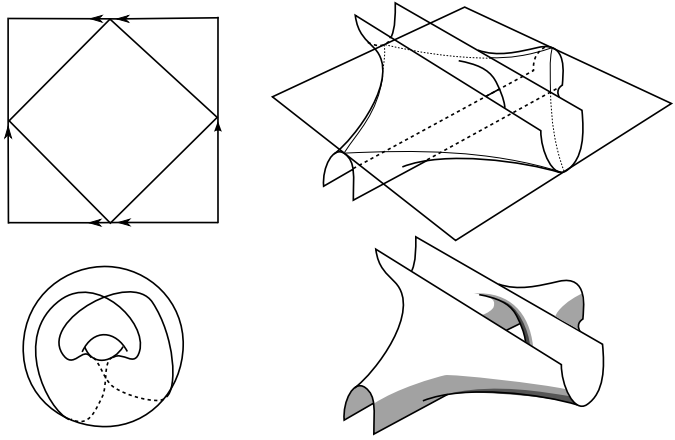}
  \caption{Extending a level surface of $f_{t_0}$ locally from a critical level in $S$ with two saddles produces a sphere with four punctures.}
  \label{torusfig}
  \end{center}
\end{figure}

Let $R$ be a regular neighborhood of $S$ and $F$ the level surface of $f_{t_0}$ that contains the critical level.  The surface $F$ intersects $R$ as shown on the right of Figure~\ref{torusfig}.  Because of the identifications at the edges of the square, this intersection is a sphere with four punctures, which we will call $U$, and a (possibly empty) collection of annuli.  The boundary loops of $U$ (and thus the boundary loops of $F \setminus R$) determine slopes in $S$ that intersect in one point.  Note that $F$ is isotopic to $\Sigma$ because it is a level surface of a sweep-out for $(\Sigma, H_1, H_2)$.

\begin{Claim}
The intersection $F \cap S$ consists the graph $U \cap S$ and a (possibly empty) collection of loops that are trivial in both $S$ and $F$.
\end{Claim}

\begin{proof}[Proof of Claim]
Let $g$ be the restriction of the sweep-out $f_{t_0}$ to the surface $S$.  Each level set of $g$ is the intersection of $S$ with a level surface of $f_{t_0}$.  There is a canonical way (up to isotopy) to identify this level surface with $\Sigma$, so each loop component of each level set of $g$ determines an isotopy class of simple closed curves in $\Sigma$.  At a central singularity in $g$, a loop corresponding to a trivial loop in $\Sigma$ is added or removed.  At a level where there is a single saddle singularity in $g$, one loop is turned into two, or vice versa by a band summing operation.  

For $t$ near $0$, these simple closed curves bound disks in $H_1$ and near $1$, they bound disks in $H_2$.  For any regular level of $g$, consisting of a number of simple closed curves, the corresponding isotopy classes of loops in $\Sigma$ are pairwise disjoint.  Because $\Sigma$ is strongly irreducible, a fixed level set of $g$ cannot determine essential loops in $\Sigma$ bounding disks on both sides.  Every regular level of $g$ contains a trivial loop in $S$, so one of the loops in $\Sigma$ determined by this regular level bounds a disk in $H_1$ or $H_2$.

The disks cannot switch from one side to the other at a critical level with a single saddle because the loops in $\Sigma$ before and after the band summing are disjoint.  Thus the switch occurs at the single level of $g$ with two saddle singularities.  In particular, the loops that limit onto this level set bound disks in opposite handlebodies (Though these disks are not disjoint.), so the remaining loops of intersection must be trivial in $F$.
\end{proof}

To complete the proof, isotope $F$ so as to remove any loops that are trivial in both $F$ and $S$.  Bachman and Derby-Talbot~\cite{bdt} pointed out that after these trivial loops are removed, $S \setminus F$ is a pair of compressing disks for $F$ whose boundaries, when made transverse, intersect in four points.  These compressing disks are on opposite sides of $F$ and are contained in the regular neighborhood $R$.  Any compressing disk for $F \setminus R$ is disjoint from each of the disks in $S \setminus F$.  Because $(\Sigma, H_1, H_2)$ is strongly irreducible and $F$ is isotopic to $\Sigma$, $\Sigma \setminus R$ must be incompressible in $M \setminus R$.  

The manifold $N \setminus R$ is homeomorphic to a pair of pants cross $S^1$.  Any incompressible surface in this manifold is either a vertical torus or a horizontal pair of pants.  The surface $F \cap (N \setminus R)$ has boundary so it is not a vertical torus and must consist of some number of horizontal pairs of pants.  The pairs of pants intersect $R$ in the loops $\partial U$.  As noted above, these loops have slopes on the boundary of the closure of $R$ that, when projected into $S$, intersect in one point.  

The first homology group of $N \setminus R$ is isomorphic to $\mathbf{Z} \times \mathbf{Z}^2$ where the first $\mathbf{Z}$ is generated by the $S^1$ factor of the pair of pants cross a circle.  The three boundary loops of a component of $F \cap (N \setminus R)$ bound a pair of pants so the sum of the homology elements they generate is zero.  The first coordinates of the two loops in $\partial R$ differ by exactly one (since they intersect in a single point in $S$) so first coordinate of the third loop must be $-1$ or $1$.    In other words, the third cuff of each pair of pants must have slope 1 or -1 in $\partial N$.

The surface $F \cap N$ is the union of a four times punctured sphere $F \cap R$ and two pairs of pants $F \cap (N \setminus R)$ so $F \cap N$ is a twice punctured genus two surface.  Thus $F \setminus N$ is an incompressible, twice punctured genus $g - 3$ surface whose boundary has slope 1 or -1 in $\partial N$.
\end{proof}

\section{Seifert fibered spaces}
\label{sfsect}

Let $M$ be a Seifert fibered space and let $c \subset M$ be a critical fiber.  The complement in $M$ of a regular neighborhood $U$ of $c$ is a surface bundle. Let $F$ be a leaf of this bundle and assume that $\partial F$ consists of a single loop in $\partial U$ that is a longitude of the solid torus $U$.  Let $\Sigma$ be union of two disjoint leaves parallel to $F$ and an annulus in $U$ connecting the boundaries of these leaves.  Each component of the complement of $\Sigma$ is homeomorphic to the union of $(\partial F) \times (0,1)$ and a regular neighborhood of $\partial F \times (0,1)$.  

Because $F$ is a surface with boundary, this set is a handlebody so $\Sigma$ determines a Heegaard splitting $(\Sigma, H_1, H_2)$.  A Heegaard splitting constructed in this way is called a \textit{horizontal Heegaard splitting}.  Recall the collection $\{T_i\}$ of vertical tori defined in Section~\ref{introsect}.  The slope that $(\Sigma, H_1, H_2)$ determines on each $T_i$ is precisely the slope of intersection between $F$ and $T_i$.  

Sedgwick~\cite{sedg} showed that a horizontal Heegaard splitting is irreducible if and only if the winding number of $c$ is greater than the least common multiple of the winding numbers of the other critical fibers.  In particular, if $\Sigma$ is strongly irreducible then the winding number of $c$ must be the largest over all the critical fibers in $M$.  Thus if $(\Sigma', H'_1, H'_2)$ is a second strongly irreducible horizontal Heegaard splitting of $M$ then $\Sigma'$ is constructed starting from the same fiber $c$.  The incompressible surface $F$ is uniquely determined (up to isotopy) by the slopes of intersection between $F$ and each $T_i$.  Thus if $\Sigma$ and $\Sigma'$ determine the same slope in each $T_i$ then they were constructed from the same $c$ and $F$, and are therefore isotopic.

\begin{proof}[Proof of Lemma~\ref{oneslopelem}]
The discussion above shows that if two strongly irreducible, horizontal Heegaard splittings determine the same slope with each $T_i$ then they are isotopic.  We will prove the converse.  Without loss of generality, assume $i$ is odd, so that $T_i \cap T_{i+1}$ is a single simple closed curve.

A regular neighborhood $N$ of $T_i \cup T_{i+1}$ is homeomorphic to a punctured torus cross an interval.  Because $T_i$ and $T_{i+1}$ is each isotopic to a union of regular fibers in $M$, we can assume that $N$ is also a union of regular fibers.  The complement in $M$ of $N$ is a Seifert fibered space so every incompressible surface in $M \setminus N$ is either a vertical torus or a horizontal incompressible surface.  The only of these surfaces that has boundary in $\partial N$ is a horizontal surface.  

Assume for contradiction $\Sigma$ determines more than one slope in $T_i$.  Then by Lemma~\ref{slopechangelem} there is an incompressible surface $F$ in the complement of $N$ that intersects the boundary in two parallel loops with slope $\pm 1$.  A horizontal incompressible surface in a Seifert fibered space is non-separating so $F$ (which is separating) must be a union of two horizontal surfaces.  The complement $M \setminus N$ is a Seifert fibered space whose fibers in $\partial M$ match the fibers in $N$, and thus have slope $\infty$ in $\partial N$.  

Each boundary component of $F$ has slope $\pm 1$ so each regular fiber of the fibrations intersects each component of $F$ in a single point.  The number of intersections of a singlular fiber with a horizontal surface is a proper integral fraction of the number of intersections with the nearby regular fibers, so $M \setminus N$ contains no singular fibers.  By construction, $N$ containes no singular fibers, so $M$ must be a circle bundle.  By construction, circle bundles do not admit strongly irreducible horizontal Heegaard splittings (see also Weidmann's classification of Heegaard splittings of circle bundles in the Appendix of~\cite{bdt}) so this contradicts the assumption that $M$ admits a strongly irreducible, horizontal Heegaard splitting.
\end{proof}

\begin{proof}[Proof of Theorem~\ref{mainthm}]
Lemma~\ref{oneslopelem} implies that a strongly irreducible horizontal Heegaard splitting is uniquely determined by the $2g$-tuple of slopes it determines with the incompressible tori $\{T_i\}$.  To show a one-to-one correspondence to $\mathbf{Z}^{2g}$ we need only show that if $M$ admits a strongly irreducible horizontal Heegaard splitting then for any $2g$-tuple there is a strongly irreducible horizontal Heegaard splittings that determines this $2g$-tuple of slopes.

If $M$ has a strongly irreducible, horizontal Heegaard splitting $(\Sigma, H_1, H_2)$, then $\Sigma$ was constructed from some critical fiber $c$ and an incompressible surface $F$ in the complement of $c$.  The critical fiber $c$ can always be taken disjoint from each $T_i$.  

Given positive integers $n$ and $i \leq g$, consider $n$ parallel copies $T_{2i}$.  The surface $F$ intersects each copy $T^j_{2i}$ of $T_{2i}$ in two simple closed curves.  Let $U$ be a regular neighborhood of a component of $F \cap T^j_{2i}$.  The intersection of $F \cup T^j_{2i}$ with $U$ is the union of a pair of annuli that intersect in a common essential loop.  There are two ways to replace these two intersecting annuli with two disjoint annuli.  If we make this replacement in the same way in each neighborhood, the resulting surface will have slope either $n$ or $-n$ in $T_{2i+1}$.  For every other $T_j$, the slopes of $F \cap T_j$ and $F' \cap T_j$ agree.  (This operation is called a \textit{Haken sum}.)  We will say that the surface with slope $n$ is the result of \textit{spinning} $F$ around $T_{2i}$ $n$ times.

Similarly, spinning $F$ around $T_{2i+1}$ changes its slope with $T_{2i}$ but not with the other vertical tori.  Thus by spinning $F$ around the vertical tori, one can construct a horizontal surface $F'$ that intersects the vertical tori $\{T_i\}$ in any $2g$-tuple.  This $F'$ has the same boundary as $F$ in $\partial N(c)$ so $F' \cup A$ is a horizontal Heegaard surface $\Sigma'$ for $M$.  There are two ways to see that $\Sigma'$ is a strongly irreducible, horizontal Heegaard surface.  First, the reader can check that there is a homeomorphism from $M$ to itself taking $\Sigma$ onto $\Sigma'$.  Second, both Heegaard splittings are constructed from the same critical fiber in $M$, so by Sedgwick's results~\cite{sedg}, both are strongly irreducible.  The Heegaard surface $\Sigma'$ determines the same $2g$-tuple of integers as $F'$, so for each $2g$-tuple of integers, there is a strongly irreducible, horizontal Heegaard splittings whose slopes in $\{T_i\}$ realize those values.
\end{proof}

\section{Double primitive knots}
\label{dpsect}

In this section we will construct a family of 3-manifolds with infinitely many non-isotopic genus three Heegaard splittings.  In the next section, we will show that for certain Seifert fibered spaces over the torus with two critical fibers, the Heegaard splittings all determine the same Nielsen classes of generators for the fundamental group.

Let $X$ be a compact, closed, orientable, irreducible 3-manifold with a genus two Heegaard splitting $(\Sigma', H'_1, H'_2)$ and let $\ell \subset \Sigma'$ be a simple closed curve such that for $i = 1,2$, $\ell$ intersects some essential, properly embedded disk $D_i \subset H'_i$ in a single point.  Such a loop is called \textit{double primitive}.  (A knot in $S^3$ that is isotopic to a double primitive loop in a genus two Heegaard splitting is called a \textit{Berge knot}.)

Let $Y \subset X$ be a regular neighborhood of a double primitive loop $\ell \subset \Sigma'$.  The intersection $\partial Y \cap \Sigma'$ is a pair of loops $\ell'_1$ and $\ell'_2$ in the torus $\partial Y$.  Define $N = T \times S^1$ where $T$ is a once punctured torus and let $x_1$,$x_2$ be points in $S^1$.  Let $M$ be the result of gluing $X \setminus Y$ to $N$ by a map that sends $\ell'_1$ to $\partial T \times \{x_1\}$ and sends $\ell'_2$ to $\partial T \times \{x_2\}$.

\begin{Lem}
\label{itsheegardlem}
The surface $\Sigma = \Sigma' \cup (T \times \{x_1,x_2\})$ is a genus three Heegaard surface for $M$.
\end{Lem}

\begin{proof}
The complement in $N$ of $(T \times \{x_1,x_2\})$ consists of two components, each of whose closure is a genus two handlebody $T \times I$ for $I \subset S^1$ one of the two intervals with endpoints $x_1, x_2$.  Let $N_1$, $N_2$ be these handlebodies.  Then each of $N_1 \cap \partial N$ and $N_2 \cap \partial N$ is an annulus.  

In $M$, the complement $M_1 = H_1 \setminus Y$ is a handlebody.  Because $\ell$ is double primitive, the intersection of $M_1$ with the closure of $Y$ is an annular neighborhood $A$ of a loop in $\partial M_1$ that intersects some properly embedded, essential disk $D \subset M_1$ in a single point.  A closed regular neighborhood $U$ in $M_1$ of $A \cup D$ is a solid torus such that $U$ intersects the closure of $M_1 \setminus U$ in a disk.  In $M$, the set $U$ is a regular neighborhood of an annulus in the boundary of $N_1$.  Thus $N_1 \cup U$ is a handlebody.  

The set $M_1 \cup N_1$ is the union of the closure of $M_1 \setminus U$ (which is a handlebody) and the handlebody $N_1 \cup U$.  The two handlebodies intersect in a disk so their union is a handlebody $H_1$.  A similar argument for $M_2$ implies that $N_2 \cup M_2$ is a handlebody $H_2$.  Thus $\Sigma = \partial H_1 = \partial H_2$ is a Heegaard surface for $M$.
\end{proof}

The Heegaard surface $\Sigma = \Sigma' \cup (T \times \{x_1,x_2\})$ determines a Heegaard splitting $(\Sigma, H_1, H_2)$ such that $H'_1 \subset H_1$ and $H'_2 \subset H_2$.  Note that the Lemma requires that for $s \in S^1$, $\partial T \times \{s\}$ is sent to the same slope in $\partial Y$ as $\ell'_1$, but there is no requirement for the slope that a loop $\{t\} \times S^1$ ($t \in \partial T$) is glued to.  Thus there are infinitely many gluings that will produce a manifold with a genus three Heegaard splitting.

\begin{Lem}
\label{itssilem}
If $(\Sigma, H_1, H_2)$ is weakly reducible then $X \setminus Y$ is a solid torus.
\end{Lem}

\begin{proof}
Because $M$ has Heegaard genus three, it cannot be a connect sum of $T^3$ with a non-trivial manifold.  If $\partial N$ is compressible then it compresses down to a sphere, which must bound a ball.  Thus if $\partial N$ is compressible then $M \setminus Y$ is a solid torus.  We will therefore assume that $\partial N$ is incompressible.

Assume further that $\Sigma$ is weakly reducible.  By Casson and Gordon's Theorem~\cite{cass:red}, if $(\Sigma, H_1, H_2)$ is weakly reducible then $\Sigma$ is reducible or $\Sigma$ compresses to a separating incompressible surface $S$ in $M$.  In the second case, each component of the complement of $S$ has a Heegaard splitting that comes from compressing $\Sigma$.  

The 3-manifold $M$ does not admit a genus two Heegaard splitting because by the main Theorem of~\cite{kob}, if a closed 3-manifold $M$ contains a separating incompressible torus and a genus two Heegaard splitting then each piece of the complement is a Seifert fibered space over a disk, an annulus or a Mobius band, the complement of a $(1,1)$ knot in a lens space or the complement of a two-bridge knot in $S^3$.  (In fact, the theorem is much stronger than this, but that's all we need.)  The component $N$ is not one of these three types so $M$ does not admit a genus two Heegaard splitting.  Because $\Sigma$ is not reducible, the weak reduction must determine a separating incompressible surface $S \subset M$.

Because $S$ is the result of compressing the genus three surface $\Sigma$ at least twice, $S$ must consist of one, two or three tori.  Because each component of $M \setminus S$ has a Heegaard splitting that comes from compressing the genus three surface $\Sigma$, each component of $M \setminus S$ has Heegaard genus at most two.  Any submanifold of $M$ containing $N$ has Heegaard genus at least three (for the same reason that $M$ has Heegaard genus at least three) so $S$ must intersect $N$.  

Any incompressible surface in $N$ is either a vertical torus or a horizontal once punctured torus.  If $S \cap N$ contains a horizontal punctured torus then $S \setminus N$ contains a disk so $\partial N$ is compressible into $X$.  This contradicts the assumption that $\partial N$ is incompressible.

Thus $S$ consists of vertical tori in $N$. Because it is separating, $S$ must consist of a union $U \subset N$ of two parallel vertical tori (each of which is non-separating).  Each component of $M \setminus U$ has a Heegaard splitting induced from $\Sigma$ and such a splitting has genus at most two.  One component is homeomorphic to a torus cross an interval.  The other component is the union of $M \setminus Y$ and a pair of pants cross an interval.

Let $Z$ be this second component.  Note that if $M \setminus Y$ is not a solid torus then the fundamental group of $M \setminus Y$ has rank at least two.  The fundamental group a pair of pants cross an interval is the direct product of $\mathbf{Z}$ and a free group $\mathbf{F}^2$ on two generators.  By Van Kampen's Theorem, the fundamental group of $Z$ is the quotient of the free product $\pi_1(M \setminus Y) * (\mathbf{F}^2 \times \mathbf{Z})$ by two relations, one of which equates an element of $\mathbf{F}^2$ to an element of $\pi_1(M \setminus Y)$ and the other equates a generator of $\mathbf{Z}$ to an element of $\pi_1(M \setminus Y)$.  There is thus a homomorphism from $\pi_1(Z)$ onto the direct product $\pi_1(M \setminus Y) \times \mathbf{Z}$.  If $Z$ admits a genus two Heegaard splitting then $\pi_1(Z)$ has rank at most 2 so $\pi_1(M \setminus Y)$ has rank at most one, implying $M \setminus Y$ is a solid torus.
\end{proof}

Let $\alpha$ and $\beta$ be essential simple closed curves in $T$ whose intersection is a single point.  Define $S_\alpha = \alpha \times S^1$ and $S_\beta = \beta \times S^1$.  Because $\Sigma$ contains $T \times x_1$ and $T \times x_2$, it determines the slope $0$ in both $S_\alpha$ and $S_\beta$.  Let $\Sigma_i$ be the result of spinning $\Sigma$ $i$ times around $S_\beta$ as in Section~\ref{sfsect}.  If two such surfaces, $\Sigma_i$ and $\Sigma_j$ ($i \neq j$), are isotopic then $\Sigma_i$ determines both the slope $i$ and the slope $j$.  By Lemma~\ref{slopechangelem}, this implies that $\Sigma$ can be isotoped to intersect $X \setminus Y$ in an incompressible, twice punctured, genus zero surface, i.e. an annulus.

Assume that $\Sigma$ can be isotoped to intersect $X \setminus Y$ in an incompressible, properly embedded annulus $A$.  If $A$ is boundary compressible then because $X$ is irreducible, $A$ must be boundary parallel.  Isotoping $A$ out of $X \setminus Y$ makes $\Sigma$ disjoint from $\partial(X \setminus Y)$.  This implies that $\partial(X \setminus Y)$ must be compressible because a Heegaard surface cannot be made disjoint from a closed incompressible surface.  Thus we have proved the following:

\begin{Lem}
\label{noannlem}
If $\partial(X \setminus Y)$ is incompressible in $M$ and $X \setminus Y$ does not contain a properly embedded, essential (i.e. incompressible and boundary incompressible) annulus then $M$ admits infinitely many non-isotopic Heegaard splittings.
\end{Lem}

Note that if $M$ is the complement of a knot $K$ in $S^3$ then Lemma~\ref{noannlem} holds whenever $K$ is not a torus knot or a cable knot.

\section{Nielsen equivalence}
\label{nielsect}

Let $(O, B_1, B_2)$ be a genus one Heegaard splitting of $S^3$ and $K \subset O$ a simple closed curve that does not bound a disk in $S^3$.  (The curve $K$ is a non-trivial torus knot.)  We can give each solid torus $B_1$, $B_2$ a Seifert fibration such that $K \subset \partial B_i$ is a fiber.  This defines a Seifert fibration of $S^3$ such that there is a regular neighborhood $Y$ of $K$ consisting of a union of fibers.  Thus the complement in $S^3$ of $Y$ is a Seifert fibered space over the disk with two singular fibers.  A regular fiber in $\partial Y$ determines the same slope as $O \cap \partial Y$.

Let $m$ be the boundary of a small disk $D$ that intersects $K$ in a single point and $O$ in a single arc.  Let $U$ be an open regular neighborhood of $m$.  Define $E$ to be the union of the twice punctured sphere $O \setminus U$ and the annulus $\partial \bar U \cap B_1$. (Here $\bar U$ is the closure of the open set $U$.)  Define $E'$ to be the union of $O \setminus U$ and $\partial \bar U \cap B_2$.

Because $m$ bounds a disk that intersects $K$ in a single point, there is a homeomorphism from $S^3$ to the result of 1 Dehn surgery on $m$ that takes $K$ onto itself.  Let $F$ be the image in this homeomorphism of $E$ and let $F'$ be the image of $E'$.  In other words, $F$ and $F'$ are the result of ``twisting'' $E$ and $E'$, respectively, about the meridian $m$.  The differences $E' \setminus E$ and $E \setminus E'$ are annuli whose union bounds the solid torus $\bar U$. The annuli meet in meridians of the solid torus.  After the Dehn surgery, $F' \setminus F$ and $F \setminus F'$ again bound the solid torus, but this time they meet along a longitude of the solid torus.  Thus there is an isotopy from $F' \setminus F$ to $F \setminus F'$.  Extending this isotopy to all of $F'$ takes $F'$ onto $F$ so $F$ and $F'$ are isotopic in $S^3$ fixing $F' \cap Y = F \cap Y$.

The surface $F$ is the result of adding a trivial handle to the genus one Heegaard surface $O$. Thus $F$ defines a genus two Heegaard splitting $(\Sigma',H'_1,H'_2)$ for $S^3$ where $\Sigma' = F$, $H'_1 = B_1 \setminus U$ and let $H'_2 = B_2 \cup \bar U$.  The intersection $D \cap H''_2 \cap B_1$ is an essential disk properly embedded in $H'_2$ and intersects $K$ in a single point.  Thus $K$ is primitive in $H'_2$.  Because $F'$ is isotopic to $F$, a similar argument for $F'$ implies that $K$ is also primitive in $H'_1$.  (The reader can check that $F$ results from taking the standard unknotting graph consisting of a core for $B_1$ and a short arc to $K$, then pushing $K$ into the resulting Heegaard surface in a way that makes it double primitive.)

Let $T$ be a once punctured torus, let $s_1$ and $s_2$ be points in $S^1$ and let $t_1$ and $t_2$ be points in $\partial T$.  Let $M$ be the result of gluing $S^3 \setminus Y$ to $T \times S^1$ so that the loops $O \cap \partial Y$ are sent to $t_1 \times S^1$ and $t_2 \times S^1$ while the loops $\Sigma' \cap \partial Y$ are sent to $\partial T \times s_1$ and $\partial T \times s_2$.  (This is possible because each component of $O \cap \partial Y$ intersects each component of $\Sigma' \cap \partial Y$ in a single point.  

Since $K$ is double primitive in $\Sigma$, Lemma~\ref{itsheegardlem} implies that the union $\Sigma = (\Sigma \setminus Y) \cup (T \times (s_1 \cup s_2))$ is the surface in a genus three Heegaard splitting $(\Sigma, H_1, H_2)$ for $M$.  By  Lemma~\ref{itssilem}, this Heegaard splitting is strongly irreducible.  Moreover, because loops of the Seifert fibration in $\partial (S^3 \setminus Y)$ are glued to vertical loops in $T \times S^1$ (which can be thought of as loops of a Seifert fibration for $N$), $M$ is a Seifert fibered space. 

Let $\alpha$ and $\beta$ be simple closed curves in $T$ that intersect in a single point.  As in the previous section, we can spin $\Sigma$ around the vertical torus $S_\beta = \beta \times S^1$ to construct an infinite family of Heegaard splittings $\{(\Sigma^i, H^i_1, H^i_2)\}$ such that $\Sigma^i$ determines the slope $i$ in the vertical torus $S_\alpha = \alpha \times S^1$.

\begin{Lem}
The Heegaard splittings $(\Sigma^i,H^i_1,H^i_2)$ and $(\Sigma^j,H^j_1,H^j_2)$ are isotopic if and only if $i = j$.  However, the generating set for $\pi_1(M)$ defined by the inclusion map $\pi_1(H^i_1) \rightarrow \pi_1(M)$ is  Nielsen equivalent to that defined by $\pi_1(H^j_1) \rightarrow \pi_1(M)$ for all $i,j$.  The generating set determined by $\pi_1(H^i_2) \rightarrow \pi_1(M)$ is Nielsen equivalent that defined by $\pi_1(H^j_2) \rightarrow \pi_1(M)$ as well.
\end{Lem}

\begin{proof}
Note that there is an essential annulus properly embedded in $M$ so Lemma~\ref{noannlem} is not enough to distinguish Heegaard splittings by their slopes.  However, this annulus intersects $\partial Y$ in the same slope as $O \cap \partial Y$ which determines the slope $\infty$ in $\partial T \times S^1$.  Thus Lemma~\ref{oneslopelem} implies that $\Sigma^i$ determines a unique slope on any vertical torus in $T \times S^1$.  Since $\Sigma_i$ determines slope $i$ and $\Sigma_j$ determines slope $j$, we conclude $\Sigma^i$ and $\Sigma^j$ are isotopic if and only if $i = j$.  

All that remains is to show that the Nielsen classes of the generators for $\pi_1(M)$ determined by these Heegaard splittings are all equivalent.  We will show this for $\Sigma_0$ and $\Sigma_1$.  A similar argument works for any $\Sigma^i$, $\Sigma^{i+1}$ and the general result follows by induction on $|i|$.

We will choose as the base point for $\pi_1(M)$ a point $p = (a,b) \in \partial T \times S^1$.  The fundamental group of the punctured torus $T \times \{b\}$ (with base point $b \in \partial T$) is a free group on two generators.  Let $x$ and $y$ be the inclusion into $\pi_1(M)$ of these generators.  We can choose $x$ and $y$ so that an arc representing $x$ intersects $S_\beta$ in a single point and is disjoint from $S_\alpha$.  Similarly, we can assume that an arc representing $y$ intersects $S_\alpha$ in a single point and is disjoint from $S_\beta$.

Let $z$ be the element of $\pi_1(M)$ defined by the loop $a \times S^1$.  Let $t$ be the element of $\pi_1(M)$ defined by a path that follows a short arc into $B_1 \subset S^3$, then follows a core of $B_1$ disjoint from the disk $D$, then follows the short arc back to $p$.  Because $z$ is determined by a regular fiber and $t$ is determined by a singular fiber of order $t$, we have $z = t^p$ for some integer $p$.

The fundamental group of $H^0_1$ is generated by $x$, $y$ and $t$, so it induces the Nielsen class $[x,y,t]$ for $\pi_1(M)$.  The only generator for $H^0_1$ that intersects $S_\beta$ is $x$.  Spinning $\Sigma^0$ around $S_\beta$ replaces $x$ with $xz = zx$ or $xz^{-1} = z^{-1}x$, while fixing $y$ and $t$.  Without loss of generality, we will assume it replaces $x$ with $xt$.  Thus $H^1_1$ determines the Nielsen class $[xz,y,t]$.  We noted above that $z = t^p$, so the new generating set is in fact $[xt^p,y,t]$, which is Nielsen equivalent to $[x,y,t]$, the generating set for $H^0_1$.  The generating sets induced by $\pi_1(H^0_1)$ and $\pi_1(H^1_1)$ (and by induction of any $\pi_1(H^i_1)$) are Nielsen equivalent.

Above, we constructed $\Sigma^0$ from the surface $F$ in the knot complement $M \setminus Y$.  Switching the roles of $B_1$ and $B_2$ in this construction switches $F$ and $F'$, so the resulting Heegaard splitting would be constructed from $F'$.  However, we noted that $F'$ is isotopic to $F$ in $M \setminus Y$.  Thus the Heegaard splitting that results from switching the roles of $B_1$ and $B_2$ is isotopic to $(\Sigma^0, H^0_2, H^0_1)$ (i.e. the same Heegaard surface, but with the order of the handlebodies switched.)  We can thus apply the argument above to $H^0_2$ and $H^1_2$, implying that the generating sets induced by $\pi_1(H^i_2)$ and $\pi_1(H^j_2)$ are Nielsen equivalent for all $i,j$.
\end{proof}

\bibliographystyle{amsplain}
\bibliography{horizontals}

\end{document}